\documentclass[12pt,a4paper]{article}
\usepackage{amssymb,amsmath}
\usepackage{amsfonts}
\usepackage{amsthm}
\usepackage{a4}
\usepackage{tikz}
\usepackage{verbatim}
\usepackage{caption}
\usepackage{float}
\newcommand{\para}{\par\vspace{.25cm}}
\newtheorem{defi}{Definition}
\newtheorem{prop}{Proposition}
\newtheorem*{theorem*}{Theorem}
\newtheorem{theorem}{Theorem}
\newtheorem{lemma}{Lemma}
\newtheorem{cor}{Corollary}
\newtheorem{remark}{Remark}

\newcommand{\Q}{\mathbb{Q}}

\newcommand{\C}{\operatorname{Cen}}
\newcommand{\I}{\operatorname{Irr}}

\begin{document}
\baselineskip 18pt \title{{\bf \ Connecting monomiality questions with the structure of rational group algebras}}
\author{ Gurmeet K. Bakshi   {\footnote {Research supported by Science and Engineering Research Board (SERB), DST, Govt. of India under the scheme Mathematical Research Impact Centric Support (sanction order no MTR/2019/001342) is gratefully acknowledged.}}\\ {\em \small Centre for Advanced Study in
Mathematics,}\\
{\em \small Panjab University, Chandigarh 160014, India}\\{\em
\small email: gkbakshi@pu.ac.in}\vspace{.2cm}\\ and \vspace{.2cm}\\ Gurleen Kaur {\footnote{Corresponding author}}\\{\em \small Department of Mathematics,}\\ {\em \small Sri Guru Gobind Singh College, Chandigarh 160019, India}\\{\em
\small email: gurleenkaur992gk@gmail.com} }
\date{}
{\maketitle}
\begin{abstract} In recent times, there has been a lot of active research on monomial groups in two different directions. While group theorists are interested in the study of their normal subgroups and Hall  subgroups, the interest of group ring theorists lie in the structure of their rational group algebras due to varied applications. The purpose of this paper is to bind the two threads together. Revisiting Dade's celebrated embedding theorem which states that a finite solvable group can be embedded inside some monomial group, it is proved here that the embedding is indeed done inside some generalized strongly monomial group. The so called generalized strongly monomial groups arose in a recent work of authors while  understanding  the algebraic structure of rational group algebras.  Still unresolved monomiality questions have been correlated by proving that all the classes of monomial groups where they have been answered are generalized strongly monomial. The study  also raises some intriguing questions weaker than those asked by Dornhoff and Isaacs in their investigations. \end{abstract}
\section{Introduction}  Throughout this paper, we assume that all groups are finite and all characters are complex characters. A group is said to be a \textit{monomial group} if each of its irreducible complex character is monomial, i.e., it is induced from a linear character of some subgroup. Dade proved that every finite solvable group can be embedded inside a monomial group (see \cite{IsaacsV}, Theorem 9.7).  A direct consequence of Dade's result is that an arbitrary subgroup of a monomial group is not necessary  monomial. In 1967, Dornhoff \cite{Dornhoff} proved that normal Hall subgroups of monomial groups are monomial and asked that whether the Hall subgroups and normal subgroups of monomial groups are again monomial. For the past half century, the focus of group theorists has been to establish connections of monomial groups with their Hall subgroups and normal subgroups. In 2005, Fukushima \cite{Fukushima} settled the question regarding Hall subgroups of monomial groups in negative by constructing a family of monomial groups where it fails. Regarding normal subgroups of monomial groups, in 1973, Dade \cite{Dade} gave an example of a monomial group of even order having a non monomial normal subgroup.  In his example, the  prime 2 played a fundamental role and so the monomiality question regarding normal subgroups of odd order monomial groups was left open. Since then several researchers have been trying hard to answer this question. Parks \cite{Parks} showed that every normal subgroup of a  nilpotent-by-supersolvable monomial group of odd order is monomial.  In \cite{Gunter}, Gunter proved that every normal subgroup of  a monomial group with a Sylow tower is monomial. In \cite{ML}, Loukaki proved that every  normal subgroup of a monomial group of order $p^{a}q^{b}$ is monomial, where $p$ and $q$ are odd primes. Later, generalizing the ideas contained in Parks's and Loukaki's work,  the theory of linear limits was introduced by Dade and Loukaki \cite{DL}. Using this theory, Chang, Zheng and Jin \cite{CZJ} generalized the work of Gunter \cite{Gunter} for a class of solvable groups with technical conditions involved in it. Extending the work of Parks \cite{Parks}, the monomiality of normal subgroups and Hall subgroups of a class of solvable groups with some restrictions on its supersolvable residual was given by Zheng and Jin \cite{ZJ}.  In \cite{IsaacsIII}, Isaacs has pointed that the problem being faced while working with monomial character is that every character inducing it need not necessary be monomial (also see Chapter 9 of \cite{IsaacsV}), and so he defined supermonomial characters. Isaacs, in \cite{IsaacsIII}, called an irreducible complex character $\chi$ of a group $G$ to be supermonomial if every character inducing $\chi$ is monomial and conjectured that every monomial group of odd order is a \textit{supermonomial group}, i.e., each of its irreducible complex character is supermonomial. Recently  Lewis, in \cite{Lewis}, proved that  if this conjecture is true, then the question about the monomiality of normal  subgroups of  odd order monomial groups  can be settled in affirmative. \par We now draw attention to the other thread on monomial groups which is connected with the study of their rational group algebras. One of the fundamental problem in group rings is to understand  the precise Wedderburn decomposition of the rational group algebra of a finite group. In this connection, Olivieri, del R{\'{\i}}o and Sim{\'o}n \cite{OdRS04} defined strongly monomial groups (including abelian-by-supersolvable groups) and gave precise description  of the simple components of their rational group algebras. Later, in \cite{BK2}, we  generalized  this concept and accordingly defined generalized strongly monomial groups. We had seen, in \cite{BK2}, that many interesting classes of monomial groups turn out to be generalized strongly monomial covering the class $\mathcal{C}$ of finite groups, introduced by Huppert (see \cite{BH}, Chapter 24), that consists of all finite groups whose each subquotient is either abelian or has a non central abelian normal subgroup. A significant result proved in \cite{BK2} is the explicit description of the primitive central idempotents and the corresponding simple components of the rational group algebra of a generalized strongly monomial group from its  subgroup structure.  For an extensive survey on this topic, we refer to \cite{B}and a quick recollection is done in section 2. \par The objective of this paper is to provide a coherent unification of these two threads. In section 3, we have revisited Dade's embedding theorem and proved the following in Theorem \ref{t1}:
\begin{quote} {\it Every finite solvable group is isomorphic to a subgroup of some generalized strongly monomial group. }
\end{quote} This motivates us to investigate the closeness of monomial groups to generalized strongly monomial groups. In sections 4-7, we examined all the groups where monomiality questions have been answered  and surprisingly found that all of them are indeed generalized strongly monomial.  To be  precise, we have shown in series of theorems (Theorems 2-10) that all the classes of groups discussed in the work of Parks \cite{Parks}, Gunter \cite{Gunter},  Loukaki \cite{ML}, Lewis \cite{Lewis}, Chang, Zheng and Jin \cite{CZJ} and Zheng and Jin \cite{ZJ} turn out to be generalized strongly monomial. Below is the list of  some well known classes of groups which have been proved to be generalized strongly monomial: \begin{itemize}  \it {
\item supermonomial groups (Theorem \ref{t3});\item  monomial groups of order $p^aq^b$ ($p$ and $q$ odd primes) (Theorem \ref{t6});  \item  monomial groups of odd order which are nilpotent-by-supersolvable (Theorem \ref{c0}); \item  monomial groups of odd order which are nilpotent-by nilpotent-by-Sylow abelian (Theorem \ref{t8})}; \item  monomial groups with Sylow towers (Theorem \ref{t9}, Corollary \ref{c4}).
\end{itemize} All of these results correlate generalized strongly monomial groups with the monomiality questions mentioned in the beginning. Also they accordingly arouse interest in investigating the counter example of Fukushima  where monomiality is not preserved by Hall subgroups, and that given by Dade  where monomiality is not inherited by normal subgroups. We have found that both of these examples are also generalized strongly monomial  (Theorem \ref{t15} and Remark \ref{r1}).\par The results proved in this paper give an evidence to expect that under fairly general circumstances a given monomial group is generalized strongly monomial. The results obtained  raise some questions weaker than those asked by leading group theorists which will be stated in section 8. Any attempt in answering to the questions raised in section 8, in addition to contribution in group theory, will also contribute to understanding of semisimple group algebras which is a problem of core interest in group rings.
	\section{Background on the algebraic structure of  rational group algebras and generalized strongly monomial groups} Let $\Q G$ be the rational group algebra of a finite group $G$. By Maschke's theorem (see \cite{PS}, Theorem 3.4.7), $\Q G$  is a semisimple ring, which means that  each  left ideal of $\Q G$ is its  direct summand. Consequently, the structure theorem due to Wedderburn and Artin (see \cite{PS}, Theorem 2.6.18) implies that $\Q G$  is uniquely written as $ \displaystyle{\oplus_{1 \leq i \leq k}} M_{n_{i}}(D_{i}) $, a  direct sum of matrix algebras over  finite dimensional division algebras $D_i$ over $\Q$. This is called the \textit{Wedderburn decomposition} of $\Q G$.  The summands $M_{n_{i}}(D_{i})$ in this decomposition are called the \textit{simple components} of $\Q G$. Furthermore, each simple component $M_{n_{i}}(D_{i})$ is a two sided ideal of $\Q G$ generated by a central idempotent  $e_i$ which  is primitive, i.e.,  $e_i$ has  the property that it can't further be written as a sum of two mutually orthogonal central idempotents of $\Q G$; such central idempotents of $\Q G$ are called \textit{primitive central idempotents}. The set $\{ e_i ~|~ 1 \leq i \leq k\}$ of all the primitive central idempotents  is uniquely determined by the rational group algebra $\Q G$  and is called the \textit{complete and irredundant set of primitive central idempotents} of $\Q G$. In practice, it is quite a hard problem to explicitly determine the complete and irredundant set of primitive central idempotents and the Wedderburn decomposition of a given semisimple group algebra $\Q G$. Furthermore,  an understanding of this problem is a tool to deal with several problems concerning group algebras. \par Given a group $G$, a natural approach to find the precise and explicit \linebreak Wedderburn decomposition  of $\Q G$ including the primitive central idempotents is using the character theory of $G$. The classical method to find  primitive \linebreak central idempotents of $\Q G$ begins with the computation of the primitive central idempotent $e(\chi)=\frac{\chi(1)}{|G|}{\sum _{g \in G}} \chi(g) g^{-1}$ of the complex group algebra $\mathbb{C} G$, followed by summing up all the primitive central idempotents of the form $e(\sigma \circ \chi)$ with $\sigma \in \operatorname{Gal}(\Q (\chi)/\Q)$ for $ \chi \in \I (G)$, where $\I (G)$ is the set of all the irreducible \linebreak complex characters of $G$, $\Q (\chi)$ is the field obtained by adjoining  all the values $\chi (g)$, for  $g \in G$, to $\Q $ and $\operatorname{Gal}(\Q (\chi)/\Q)$ is the Galois group of $\Q (\chi)$ over $\Q$ (see \cite{Yamada} for details). The primitive central idempotent so obtained, i.e., $\sum_{\sigma \in \operatorname{Gal}(\mathbb{Q}(\chi)/ \mathbb{Q})}e(\sigma \circ \chi) $ of $\Q G$, is commonly denoted by $e_{\mathbb{Q}}(\chi)$ and is called the \textit{primitive central idempotent of $\Q G$ realized by $\chi$}.  This approach  has computational difficulty because  it is difficult to determine $\operatorname{Gal}(\mathbb{Q}(\chi)/ \mathbb{Q})$ for a given $\chi$ and moreover the complete information about the character table of $G$ may not be known. \par For a monomial character $\chi$, Olivieri, del R{\'{\i}}o and Sim{\'o}n \cite{OdRS04} gave the expression of $e_{\Q}(\chi) $  that avoids the knowledge of $\operatorname{Gal}(\mathbb{Q}(\chi)/ \mathbb{Q})$. They found that if $\chi$ is an irreducible complex character of $G$ which is induced from a linear character $\lambda$ of some subgroup $H$ of $G$ with $K$ as its kernel, i.e., $K= \ker\lambda$, then the expression  of $e_{\Q}(\chi) $  can be written in terms of $G$, $H$ and $K$. In \cite{OdRS04}, a pair $(H,K) $ (with $H/K$ cyclic) of subgroups of $G$ with the property that  linear character of $H$ with kernel $K$  induces irreducibly to $G$ is termed as a Shoda pair of $G$. The term Shoda pair is in honor of Shoda who gave the criterion for the irreducibility of induced monomial characters (see \cite{JdR}, Corollary 3.2.3). More precisely, a Shoda pair is defined as follows:
\begin{defi}{\rm  Shoda Pair (\cite{OdRS04}, Definition 1.4)}  A Shoda pair $(H, K)$ of $G$ is a pair of subgroups of $G$ satisfying the following:\\
	{\rm(i)} $K \unlhd H,$ $ H/K$ is cyclic;\\
	{\rm(ii)} if $g \in G$ and $[H, g] \cap H \subseteq K,$ then $ g \in H.$
\end{defi}
\noindent For $K\unlhd H\leq G$, define:$$\widehat{H}:=\frac{1}{|H|}\displaystyle\sum_{h \in H}h,$$ $$\varepsilon(H,K):=\left\{\begin{array}{ll}\widehat{K}, & \hbox{$H=K$;} \\\prod(\widehat{K}-\widehat{L}), & \hbox{otherwise,}\end{array}\right.$$ where $L$ runs over all the minimal normal subgroups of $H$ containing $K$ properly, and $$e(G,H,K):= {\rm~the~sum~of~all~the~distinct~}G{\rm {\tiny{\operatorname{-}}} conjugates~of~}\varepsilon(H,K).$$ \para \noindent  If $(H, K)$ is a Shoda pair of $G$ and $\lambda$ is a linear character of $H$ with kernel $K$, then Theorem 2.1 of \cite{OdRS04} tells that $e_{\mathbb{Q}}(\lambda^G)$ is a rational multiple (unique) of $e(G, H, K)$. It may be pointed out that this knowledge was not enough to accomplish the task of determining  the structure of the corresponding simple component $\Q G e_{\Q}(\lambda^G) $. However, this task was achieved by Olivieri, del R{\'{\i}}o and Sim{\'o}n \cite{OdRS04} by imposing certain constraints on the Shoda pair $(H, K)$. The Shoda pairs with these added constraints are termed as strong Shoda pairs and are defined as follow:
\begin{defi}{\rm Strong Shoda Pair (\cite{OdRS04}, Definition 3.1, Proposition 3.3)}  A Shoda pair $(H, K)$ of $G$ is called a strong Shoda pair if \\ {\rm(i)}  $H \unlhd N_{G}(K);$ \\ {\rm(ii)}
	$\varepsilon(H, K)\varepsilon(H, K)^{g}=0$ for all $ g \in G \setminus N_{G}(K),$  where $\varepsilon(H, K)^{g}= g^{-1}\varepsilon(H, K) g$.
\end{defi}

\begin{defi}{\rm Strongly monomial character (see \cite{JdR}, p.104)} An irreducible complex character $\chi$ of a group $G$ is said to be a strongly monomial  character of $G$ if $\chi = \lambda^G$ for a linear character $\lambda$ of a subgroup $H$ with kernel $K$ such that $(H, K)$ is a strong Shoda pair of $G$.
\end{defi}
 \begin{defi}{\rm Strongly monomial group (\cite{JdR}, p.104)} A group $G$  is said to be a strongly monomial group if each of its irreducible complex character is strongly monomial.\end{defi} \noindent It is proved in Theorem  4.4 of \cite{OdRS04} that every abelian-by-supersolvable group is strongly monomial.\par Recently, in \cite{BK2}, we have given  a generalization of the concept of strong Shoda pairs, and  correspondingly strongly monomial groups, and  called them generalized strong Shoda pairs and generalized  strongly monomial groups, respectively. \begin{defi} {\rm Generalized strong Shoda pair  (\cite{BK2}, p.422):} If $(H,K)$ is a Shoda pair of $G$ and $\lambda$ is a linear character of $H$ with kernel $K$, then we say that the pair $(H,K)$ is  a \textit{generalized strong Shoda pair} of $G$ if there is a chain $H=H_{0}\leq H_{1}\leq \cdots \leq H_{n}=G$ (called strong inductive chain from $H$ to $G$) of subgroups of $G$ such that the following conditions hold for all $ 0 \leq i \leq n-1$: \begin{description} \item [(i)] $H_i \unlhd \operatorname{Cen}_{H_{i+1}}(e_{\mathbb{Q}}(\lambda^{H_i}))$; \item [(ii)] the distinct $H_{i+1}$-conjugates of $e_{\mathbb{Q}}(\lambda^{H_i})$ are mutually orthogonal.\end{description}
	\end{defi} \noindent  It is easy to see that every strong Shoda pair $(H,K)$ of $G$ is a generalized strong Shoda pair (with strong inductive chain $H \leq G$) because, by (\cite{OdRS04}, Lemma 1.2, Proposition 3.3), $e_{\Q}(\lambda) = \varepsilon(H, K)$ and $\C_{G}(\varepsilon(H, K)) = N_{G}(K)$.
\begin{defi} {\rm Generalized strongly monomial  character:} An irreducible complex character $\chi$ of a group $G$ is said to be a generalized strongly monomial character of $G$ if $\chi = \lambda^G$ for a linear character $\lambda$ of some subgroup $H$ with kernel $K$ such that $(H, K)$ is a generalized strong Shoda pair of $G$. \end{defi} \begin{defi}{\rm Generalized strongly monomial group (\cite{BK2}, p.423):} A group $G$ is said to be a \textit{generalized strongly monomial group} if each of its irreducible complex character is generalized strongly monomial.
\end{defi}\noindent Beside strongly monomial groups, a list of important families of groups which are generalized strongly monomial is produced in Theorem 1 of \cite{BK2}. The beauty of this class of generalized strongly monomial groups is due to the explicit description of the primitive central idempotents and the corresponding simple components of their rational group algebras. This work appears in \cite{BK2}.
\section{Embedding  solvable group inside generalized strongly monomial group} One of the significant property of monomial groups proved by Dade (see \cite{IsaacsV}, Theorem 9.7) is that every solvable group can be embedded inside some monomial group. Given a finite solvable group $G$, consider a subnormal series  $\{e\}= G_0 \unlhd G_1 \unlhd \cdots \unlhd G_n=G$  with factor groups $G_i/G_{i-1}$, $ 1\leq i \leq n$,  of prime order.  Denote $G_i/G_{i-1}$ by $C_i$. The crucial steps towards Dade's proof are the following: \begin{enumerate} \item  $G$ is isomorphic to a subgroup of $(((C_1 \wr C_2)\wr C_3) \wr \cdots )\wr C_n$, where $A \wr B$ denotes the wreath product of $A$ by $B$. This is a direct consequence of Lemma 9.6 of  \cite{IsaacsV}; \item  the wreath product of a monomial group by a cyclic group of prime order is monomial (see Lemma 9.5 of \cite{IsaacsV}). \end{enumerate} In this section, we will show  in Proposition \ref{p1} that the wreath product of a generalized strongly monomial group by a cyclic group of prime order is  generalized strongly monomial. Consequently, prime order groups being generalized strongly monomial, it follows that $(((C_1 \wr C_2)\wr C_3) \wr \cdots )\wr C_n$ is generalized strongly monomial and hence we obtain the  following:
\begin{theorem}\label{t1}Every finite solvable group is isomorphic to a subgroup of some \linebreak generalized strongly monomial group.
\end{theorem} \noindent As said earlier, all we  require is the following:
\begin{prop}\label{p1} If $W = A \wr C$  is the wreath product of  a generalized strongly monomial group $A$  by  a cyclic group $C$ of prime order,  then $W$ is  generalized strongly monomial.
\end{prop} \noindent To prove Proposition \ref{p1}, a bit of preparation is needed and the following lemma plays a key role.
\begin{lemma}\label{L1} Let $W=G \rtimes C$ be the semidirect product of a  finite group $G$ by a cyclic group   $C$ of prime order $p$.   Let  $\psi \in \I (G)$ be a generalized strongly monomial character which is $C$-invariant so that it extends to $W$. Assume $\psi = \lambda^G$ for some linear character  $\lambda$ of a subgroup $H$ of $G$ and that there is a strong inductive chain $ H=H_0 \leq H_1 \leq \cdots \leq H_n=G$ from $H$ to $G$ such that $C$ normalizes $H_i$ and stabilizes $\lambda^{H_i}$ for all $0 \leq i \leq n$. Then every extension of $\psi$ to $W$ is generalized strongly monomial.
\end{lemma}\noindent {\bf Proof.} Let  $\lambda^G=\psi$ has a strong inductive chain $ H=H_0 \leq H_1 \leq \cdots \leq H_n=G$  from $H$ to $G$  so that $C$ normalizes $H_i$ and stabilizes $\lambda^{H_i}$ for all $0 \leq i \leq n$. Observe that for any $i$, $0 \leq i \leq n$,  $H_iC$  is  a subgroup of $W$ and $\lambda^{H_i}$ extends to $H_i C$ and that also in $p$-ways. In particular, for $i=0,$ it gives that $\lambda$ extends to $HC$. Let $\varphi_1$, $\varphi_2$, $\cdots$, $\varphi_p$ be  all the extensions  of  $\lambda$ to $HC$. Now, for any $j$, $(\varphi_j^{H_iC})_{H_i} = ((\varphi_j)_H)^{H_i}= \lambda^{H_i}$  and so  $\varphi_j^{H_iC}$ is an extension of $\lambda^{H_i}$. Furthermore, if $\chi$ is any extension of $\lambda^{H_i}$  to $H_i C$, then by Gallaghar's theorem, $\chi = \beta \varphi_1^{H_iC}$ for some $\beta \in \I (H_i C/H_i)$. Since $\beta \varphi_1^{H_iC}= (\beta_{HC}\varphi_1)^{H_i C}$  and the restriction of  $\beta_{HC}\varphi_1$ to $H$ is $\lambda$,  $\beta_{HC}\varphi_1$ equals $\varphi_j$ for some $j$. Therefore,  $\chi = \varphi_j^{H_iC}$ and hence $\varphi_j^{H_iC}$ for $1 \leq j \leq p$ are precisely all the extensions of  $\lambda^{H_i}$ to $H_i C$. For $i=n$, it gives $\varphi_j^{W}$ for $1 \leq j \leq p$ are precisely all the extensions of  $\lambda^{G}$ to $W$. Thus, to prove the lemma, we need to show that $\varphi_j^W$, for  all $1 \leq j\leq p$, are generalized strongly monomial. \par Let $\varphi$ be any of the $\varphi_j$, $1 \leq j \leq p$. We will show that for the character $\varphi^W$, $ HC=H_0 C\leq H_1C \leq \cdots \leq H_nC=GC=W$  is a strong inductive chain from $HC$ to $W$, i.e.,  the following  hold for all $0 \leq i \leq n-1$:
\begin{description}
	\item[(i)] $ H_iC \unlhd \C_{H_{i+1}C}(e_{\mathbb{Q}}(\varphi^{H_iC}));$
\item [(ii)] distinct $H_{i+1}C$-conjugates of  $e_{\mathbb{Q}}(\varphi^{H_iC})$ are mutually orthogonal.
\end{description}
This will be proved in steps: \vspace{.2cm}\\
 \underline{\textbf{Step 1}} $ H_iC \unlhd  \C_{H_{i+1}}(e_{\mathbb{Q}}(\lambda^{H_i}))C.$ \vspace{.2cm}\\ Note that $C$ normalizes $\C_{H_{i+1}}(e_{\mathbb{Q}}(\lambda^{H_i}))$. For  if $x \in \C_{H_{i+1}}(e_{\mathbb{Q}}(\lambda^{H_i}))$ and $y \in C$, then $(\lambda^{H_{i}})^{y^{-1}xy}= (\lambda^{H_{i}})^{xy}= (\sigma \circ \lambda^{H_{i}})^{y}$ for some $\sigma \in \operatorname{Gal}(\mathbb{Q}(\lambda^{H_i})/\mathbb{Q})$ which is further equal to $\sigma \circ \lambda^{H_{i}}$, since $C$ stabilizes $\lambda^{H_{i}}$.  Thus  $e_{\mathbb{Q}}(\lambda^{H_i})^{y^{-1}xy} = e_{\mathbb{Q}}((\lambda^{H_{i}})^{y^{-1}xy}) = e_{\mathbb{Q}}(\sigma \circ \lambda^{H_{i}})= e_{\mathbb{Q}}(\lambda^{H_{i}}) ,$  implying that $y^{-1}xy \in  \C_{H_{i+1}}(e_{\mathbb{Q}}(\lambda^{H_i}))$. \par Since $C$ and $\C_{H_{i+1}}(e_{\mathbb{Q}}(\lambda^{H_i}))$ both normalizes $H_i$, all we need to show is that if $a \in C$ and $x \in \C_{H_{i+1}}(e_{\mathbb{Q}}(\lambda^{H_i}))$, then $[a, x] \in H_{i}C$. We will indeed show that $[a, x] \in H$. As $C$ normalizes $\C_{H_{i+1}}(e_{\mathbb{Q}}(\lambda^{H_i}))$, $[a, x] \in \C_{H_{i+1}}(e_{\mathbb{Q}}(\lambda^{H_i}))$. In view of  (\cite{BK2}, Theorem 2), $[a, x] $ will belong to $H$ if it stabilizes $\lambda^{H_i}$. Let us see this now. Since $H_i \unlhd  \C_{H_{i+1}}(e_{\mathbb{Q}}(\lambda^{H_i}))$ and $C~{\rm stabilizes} ~\lambda^{H_i}$, the character  $(\lambda^{H_i})^{axa^{-1}x^{-1} }$  equals to $(\lambda^{H_i})^{xa^{-1}x^{-1}}. $ Further,  $(\lambda^{H_i})^{x} = \sigma \circ \lambda^{H_i}$ for some $\sigma$ belonging to $\operatorname{Gal}(\mathbb{Q}(\lambda^{H_i})/\mathbb{Q})$, as $x \in \C_{H_{i+1}}(e_{\mathbb{Q}}(\lambda^{H_i}))$. Therefore, $(\lambda^{H_i})^{axa^{-1}x^{-1} }$=$(\lambda^{H_i})^{xa^{-1}x^{-1} }$=$(\sigma \circ \lambda^{H_i}) ^{a^{-1}x^{-1}}$= \linebreak $(\sigma \circ \lambda^{H_i}) ^{x^{-1}}$=$\sigma \circ (\sigma^{-1} \circ \lambda^{H_i})$=$\lambda^{H_i}.$ This proves step 1. \vspace{.2cm}\\  \underline{\textbf{Step 2}} If $ x \in  H_{i+1}C\setminus \C_{H_{i+1}}(e_{\mathbb{Q}}(\lambda^{H_i}))C$, then  $e_{\mathbb{Q}}(\varphi^{H_iC})e_{\mathbb{Q}}(\varphi^{H_iC})^x=0$. \vspace{.2cm}\\Consider $ x = yz$, where $y \in H_{i+1}\setminus \C_{H_{i+1}}(e_{\mathbb{Q}}(\lambda^{H_i}))$ and $ z \in C$. Then $$e_{\mathbb{Q}}(\varphi^{H_iC})e_{\mathbb{Q}}(\varphi^{H_iC})^{yz} = e_{\mathbb{Q}}(\varphi^{H_iC}) e_{\mathbb{Q}}(\varphi^{H_iC})^{zz^{-1}yz} = e_{\mathbb{Q}}(\varphi^{H_iC})e_{\mathbb{Q}}(\varphi^{H_iC})
^{z^{-1}yz}.$$ Since $C$ normalizes both  $H_{i+1}$ and  $\C_{H_{i+1}}(e_{\mathbb{Q}}(\lambda^{H_i}))$, we have $z^{-1}yz$ belongs to $ H_{i+1}$ but is not in  $\C_{H_{i+1}}(e_{\mathbb{Q}}(\lambda^{H_i})) $ and so \begin{equation}\label{e1}
e_{\mathbb{Q}}(\lambda^{H_i})e_{\mathbb{Q}}(\lambda^{H_i})^{z^{-1}yz} =0. \end{equation}  Since $\lambda^{H_i}$ extends to $\varphi^{H_iC}$, we have by (\cite{BK2}, Lemma 1) that \begin{equation}\label{e2} e_{\mathbb{Q}}(\lambda^{H_i})e_{\mathbb{Q}}(\varphi^{H_iC}) = e_{\mathbb{Q}}(\varphi^{H_iC}) = e_{\mathbb{Q}}(\varphi^{H_iC})e_{\mathbb{Q}}(\lambda^{H_i}).
\end{equation} From eqns \ref{e1} and \ref{e2}, it turns out that $e_{\mathbb{Q}}(\varphi^{H_iC})e_{\mathbb{Q}}(\varphi^{H_iC})^{z^{-1}yz}=0,$  as desired. \vspace{.2cm}\\  \underline{\textbf{Step 3}} $\C_{H_{i+1}C}(e_{\mathbb{Q}}(\varphi^{H_iC})) \leq
\C_{H_{i+1}}(e_{\mathbb{Q}}(\lambda^{H_i}))C.$  \vspace{.2cm}\\This  follows  from  step 2. \vspace{.2cm}\\ \underline{\textbf{Step 4}} $ H_iC \unlhd \C_{H_{i+1}C}(e_{\mathbb{Q}}(\varphi^{H_iC})).$ \vspace{.2cm}\\This is an immediate consequence of steps 1 and 3. \vspace{.2cm}\\ \underline{\textbf{Step 5}} $e_{\mathbb{Q}}(\varphi^{H_iC})e_{\mathbb{Q}}(\varphi^{H_iC})^x=0$ for all $ x \in \C_{H_{i+1}}(e_{\mathbb{Q}}(\lambda^{H_i}))C\setminus \C_{H_{i+1}C}(e_{\mathbb{Q}}(\varphi^{H_iC}))$. \para \noindent Let  $ x \in   \C_{H_{i+1}}(e_{\mathbb{Q}}(\lambda^{H_i}))C\setminus \C_{H_{i+1}C}(e_{\mathbb{Q}}(\varphi^{H_iC}))$. By  step 1,  $ H_iC$ is a normal subgroup of  $\C_{H_{i+1}}(e_{\mathbb{Q}}(\lambda^{H_i}))C$. Thus both  $e_{\mathbb{Q}}(\varphi^{H_iC})$ and $e_{\mathbb{Q}}(\varphi^{H_iC})^{x}$ are primitive central idempotents of $H_i C$ which must be  either same or mutually orthogonal.  Since $ x \not\in \C_{H_{i+1}C}(e_{\mathbb{Q}}(\varphi^{H_iC}))$, they can't be same and hence mutually orthogonal.  This proves step 5 and completes the proof of the lemma.\qed
\begin{lemma}\label{L2} A finite  direct product of generalized strongly monomial groups is  \linebreak generalized strongly monomial.\end{lemma} \noindent {\bf Proof.} It is enough to show that  if  $\chi_{1}$ and $\chi_{2}$ are generalized strongly monomial characters of arbitrary groups  $G_1$ and $G_2$ respectively, then   $\chi_{1} \times \chi_{2}$  is a generalized strongly monomial character of their direct product $G_{1} \times G_{2}$.  Denote $G_{1} \times G_{2}$  by $G$ and $\chi_{1} \times \chi_{2}$ by $\chi$.  Since $\chi_{1}$ is a generalized strongly monomial character of $G_{1}$, there is a strong inductive chain $H=H_{0}\leq H_{1}\leq \cdots \leq H_{t}=G_{1}$ from $H$ to $G_{1}$ and a linear character $\lambda$ of $H$ with $\lambda^{G_{1}} = \chi_{1}$, $H_{i} \unlhd \operatorname{Cen}_{H_{i+1}}(e_{\mathbb{Q}}(\lambda^{H_{i}}))$ and distinct $H_{i+1}$-conjugates of $e_{\mathbb{Q}}(\lambda^{H_{i}})$ are mutually orthogonal for all $0\leq i \leq t-1$. Similarly for $\chi_{2}$, there is a strong inductive chain $L=L_{0} \leq L_{1} \leq \cdots \leq L_{k}=G_{2}$ from $L$ to $G_{2}$ and a linear character $\vartheta$ of $L$ with $\vartheta^{G_{2}} = \chi_{2}$, $L_{i} \unlhd \operatorname{Cen}_{L_{i+1}}(e_{\mathbb{Q}}(\vartheta^{L_{i}}))$ and distinct $L_{i+1}$-conjugates of $e_{\mathbb{Q}}(\vartheta^{L_{i}})$ are mutually orthogonal for all $0\leq i \leq k-1$. We can assume that $t \leq k$. We will show that for $\chi$, $H \times L =H_{0}\times L_{0} \leq H_{1} \times L_{1} \leq \cdots \leq H_{t}\times L_{t} \leq \cdots \leq H_{k} \times L_{k} = G_{1} \times G_{2},$ where $H_{i} =G_{1}$ for all $t \leq i \leq k$ is a strong inductive chain from $H \times L$ to $G_{1}\times G_{2}$. \par  Observe that $H_{i} \times L_{i} \unlhd \operatorname{Cen}_{H_{i+1}}(e_{\mathbb{Q}}(\lambda^{H_{i}})) \times
\operatorname{Cen}_{L_{i+1}}(e_{\mathbb{Q}}(\vartheta^{L_{i}}))$ for all $0 \leq i \leq k-1$, as $H_{i} \unlhd \operatorname{Cen}_{H_{i+1}}(e_{\mathbb{Q}}(\lambda^{H_{i}}))$ and $L_{i} \unlhd \operatorname{Cen}_{L_{i+1}}(e_{\mathbb{Q}}(\vartheta^{L_{i}}))$. \par Firstly, we are going to show that  if  $(x,y) \in H_{i+1} \times L_{i+1}$ and doesn't \linebreak belong to $\operatorname{Cen}_{H_{i+1}}(e_{\mathbb{Q}}(\lambda^{H_{i}})) \times
\operatorname{Cen}_{L_{i+1}}(e_{\mathbb{Q}}(\vartheta^{L_{i}}))$, then $e_{\mathbb{Q}}(\lambda^{H_{i}} \times \vartheta^{L_{i}})^{(x,y)}$ and \linebreak $e_{\mathbb{Q}}(\lambda^{H_{i}} \times \vartheta^{L_{i}})$ are mutually orthogonal. Let $(x,y)$ belong to  $H_{i+1} \times L_{i+1}$ and doesn't belong to $\operatorname{Cen}_{H_{i+1}}(e_{\mathbb{Q}}(\lambda^{H_{i}})) \times
\operatorname{Cen}_{L_{i+1}}(e_{\mathbb{Q}}(\vartheta^{L_{i}}))$. Then  either $x \notin \operatorname{Cen}_{H_{i+1}}(e_{\mathbb{Q}}(\lambda^{H_{i}}))$ or $y \not\in \operatorname{Cen}_{L_{i+1}}(e_{\mathbb{Q}}(\vartheta^{L_{i}}))$. W.l.o.g. assume that $x \not\in \operatorname{Cen}_{H_{i+1}}(e_{\mathbb{Q}}(\lambda^{H_{i}}))$. Then \begin{equation}\label{f1}
e_{\mathbb{Q}}(\lambda^{H_{i}})^x e_{\mathbb{Q}}(\lambda^{H_{i}}) = 0.
\end{equation}Consider the restriction of  $\lambda^{H_{i}} \times \vartheta^{L_{i}}$ on $H_{i} \times \{e\}$ and observe that the restriction is homogeneous and $\lambda^{H_{i}} \times 1_{\{e\}}$ is the irreducible component, where $1_{\{e\}}$ is the principal character of the trivial subgroup $\{e\}$.  We now  use  Lemma 1 of \cite{BK2} and obtain that  \begin{equation}\label{f2}
e_{\mathbb{Q}}(\lambda^{H_{i}} \times \vartheta^{L_{i}})e_{\mathbb{Q}}(\lambda^{H_{i}} \times 1_{\{e\}})= e_{\mathbb{Q}}(\lambda^{H_{i}} \times \vartheta^{L_{i}})=e_{\mathbb{Q}}(\lambda^{H_{i}} \times 1_{\{e\}})e_{\mathbb{Q}}(\lambda^{H_{i}} \times \vartheta^{L_{i}}). \end{equation}  View the rational group algebra of $H_{i} \times L_{i}$, as the tensor product $ \mathbb{Q}H_{i} \otimes_{\mathbb{Q}}\mathbb{Q}L_{i} $ and note that  $e_{\mathbb{Q}}(\lambda^{H_{i}} \times 1_{\{e\}})= e_{\mathbb{Q}}(\lambda^{H_{i}})\otimes  e $.  Also $e_{\mathbb{Q}}(\lambda^{H_{i}} \times 1_{\{e\}})^{(x,y)}= e_{\mathbb{Q}}(\lambda^{H_{i}})^{x}\otimes  e $. Using eqn \ref{f1},  we see that  \begin{equation}\label{f3} e_{\mathbb{Q}}(\lambda^{H_{i}} \times 1_{\{e\}})^{(x,y)} e_{\mathbb{Q}}(\lambda^{H_{i}} \times 1_{\{e\}}) = e_{\mathbb{Q}}(\lambda^{H_{i}})^{x} e_{\mathbb{Q}}(\lambda^{H_{i}})\otimes  e = 0 \otimes  e =0.
\end{equation} Now, in view of eqn \ref{f2}, $ e_{\mathbb{Q}}(\lambda^{H_{i}} \times \vartheta^{L_{i}})^{(x,y)} e_{\mathbb{Q}}(\lambda^{H_{i}} \times \vartheta^{L_{i}}) $ is equal to  $$e_{\mathbb{Q}}(\lambda^{H_{i}} \times \vartheta^{L_{i}})^{(x,y)}e_{\mathbb{Q}}(\lambda^{H_{i}} \times 1_{\{e\}})^{(x,y)}e_{\mathbb{Q}}(\lambda^{H_{i}} \times 1_{\{e\}})  e_{\mathbb{Q}}(\lambda^{H_{i}} \times \vartheta^{L_{i}}), $$ which is zero, using eqn \ref{f3}, as desired.\par As a consequence of what has been shown in the above paragraph, we have
 that $\operatorname{Cen}_{H_{i+1} \times L_{i+1}}(e_{\mathbb{Q}}(\lambda^{H_{i}} \times \vartheta^{L_{i}}))$ is subgroup of $  \operatorname{Cen}_{H_{i+1}}(e_{\mathbb{Q}}(\lambda^{H_{i}})) \times
\operatorname{Cen}_{L_{i+1}}(e_{\mathbb{Q}}(\vartheta^{L_{i}}))$ and hence $H_{i} \times L_{i}$ is normal in $\operatorname{Cen}_{H_{i+1}\times L_{i+1}}(e_{\mathbb{Q}}(\lambda^{H_{i}} \times \vartheta^{L_{i}}))$. \par It now only remains to show that if $(x,y) \in
\operatorname{Cen}_{H_{i+1}}(e_{\mathbb{Q}}(\lambda^{H_{i}})) \times
\operatorname{Cen}_{L_{i+1}}(e_{\mathbb{Q}}(\vartheta^{L_{i}})) $ but $(x,y) \not\in \operatorname{Cen}_{H_{i+1}\times L_{i+1}}(e_{\mathbb{Q}}(\lambda^{H_{i}} \times \vartheta^{L_{i}}))$, then $e_{\mathbb{Q}}(\lambda^{H_{i}} \times \vartheta^{L_{i}})^{(x,y)}$ and $e_{\mathbb{Q}}(\lambda^{H_{i}} \times \vartheta^{L_{i}})$ are mutually orthogonal. Since $H_{i} \times L_{i} \unlhd \operatorname{Cen}_{H_{i+1}}(e_{\mathbb{Q}}(\lambda^{H_{i}})) \times
\operatorname{Cen}_{L_{i+1}}(e_{\mathbb{Q}}(\vartheta^{L_{i}}))$, $e_{\mathbb{Q}}(\lambda^{H_{i}} \times \vartheta^{L_{i}})$ and $e_{\mathbb{Q}}(\lambda^{H_{i}} \times \vartheta^{L_{i}})^{(x,y)}$ are primitive central idempotents of  the rational group algebra of $H_{i} \times L_{i}$, and are also distinct, as $(x,y)$ doesn't belong to $\operatorname{Cen}_{H_{i+1}\times L_{i+1}}(e_{\mathbb{Q}}(\lambda^{H_{i}} \times \vartheta^{L_{i}}))$. Therefore, they must be mutually orthogonal. This finishes the proof. \qed
\begin{lemma}\label{L3}  If  $W = A \wr C$  is the wreath product of  a generalized strongly monomial group $A$  by  a cyclic group $C$ of order $n$ and if $\chi \in \I (W)$ restricts to an irreducible character $\psi$ of the base group $A\times \cdots \times A$ ($n$ copies) of $W$ then $\psi$ is generalized strongly monomial and has a strong inductive chain satisfying the conditions stated in Lemma \ref{L1}.
\end{lemma}\noindent {\bf Proof.} Since $A$ is a generalized strongly monomial group, by Lemma \ref{L2}, we have that the base group $B$, i.e., $A \times A \times \cdots \times A$ ($n$ copies) is  generalized strongly monomial. Therefore, $\psi$ is a generalized strongly monomial character of $B$. We can write $\psi$ as  $\psi_{1} \times \psi_{2}\times \cdots \times \psi_{n},$ where $\psi_{i} \in \operatorname{Irr}(A)$ for all $1 \leq i \leq n$. The way elements of $C$ act on $B$ and $\psi$ being $C$-invariant, it turns out that all $\psi_{i}$ are equal, say equal to $\psi'$. Since $\psi'$ is a generalized strongly monomial character of $A$, we have that $\psi' = \lambda^{A}$ and  $H=H_{0} \leq H_{1} \leq \cdots \leq H_{r}=A$ is a strong inductive chain from $H$ to $A$. Now $\psi' = \lambda^{A}$ implies that  $\psi = (\lambda \times \cdots \times \lambda)^{B}$. Furthermore,  $H \times \cdots \times H \leq H_{1}\times \cdots \times H_{1} \leq \cdots \leq H_{r} \times \cdots \times H_{r}$ ($n$ copies) is a strong inductive chain of $\psi$ as shown in the proof of Lemma \ref{L2}.  From the action of $C$ on $B$, it follows that $C$ normalizes  $H_{i}\times \cdots \times H_{i}$ and stabilizes $(\lambda \times \cdots \times \lambda)^{H_{i}\times \cdots \times H_{i}}$, as desired.\qed \vspace{.5cm}\\
{\bf Proof of Proposition \ref{p1}.} If $\chi \in \I (W)$, then there is an irreducible character $\psi$ of  the base group $B= A\times \cdots \times A$ ($p$-copies) of $W$ such that either $\chi_B = \psi$ or $\chi = \psi^G$.  \par If $\chi_B =\psi$, then by Lemma \ref{L3}, $\psi$  is generalized strongly monomial and has a strong inductive chain satisfying conditions stated in Lemma \ref{L1}. Consequently Lemma \ref{L1} yields that $\chi$ is generalized strongly monomial. \par Suppose $\chi = \psi^G$. In this case, it is easy to see that for the character $\psi$ if $ H=H_0 \leq H_1 \leq \cdots \leq H_n=B$ is a strong inductive chain from $H$ to $B$, then $\chi$ has a strong inductive chain $H=H_0 \leq H_1 \leq \cdots \leq H_n=B \leq W$ from $H$ to $W$ and hence the result is proved.  \qed
\section{Supermonomial groups}
Let $\chi \in \operatorname{Irr}(G)$. Let  $N \unlhd G$ and let $\psi \in \operatorname{Irr}(N | \chi_N)$, i.e., $\psi$ is  an irreducible constituent of $\chi_{N}$.  Let $T= I_{G}(\psi),$ i.e., the inertia group of $\psi$ in $G$. By Clifford's correspondence theorem (\cite{IM}, Theorem 6.11), there is a unique Clifford correspondent  $\theta \in \operatorname{Irr}(T)$ of $\chi$ w.r.t. $\psi$ such that $\chi=\theta^{G}$ and $\theta_{N}=e\psi$, where $e=\langle \theta_{N},\psi \rangle$. Now repeat this process with $\chi$ replaced by $\theta$ and a normal subgroup $N$ of $G$ replaced by a normal subgroup $S$ of $T$ and consider the Clifford correspondent of $\theta$ w.r.t. an irreducible constituent of $\theta_{S}$. This iterative process can be continued. A character $\theta$ obtained through any number of such iterations is called a \textit{compound Clifford correspondent} of $\chi$. Those compound Clifford correspondents of $\chi$ which themselves have no proper Clifford correspondents (in other words which are quasi-primitive)  are called the \textit{stabilizer limits} of $\chi$. This terminology was introduced by Isaacs in \cite{IsaacsII}.
\begin{prop}\label{t2} Let $G$ be a solvable group,  $\chi \in \operatorname{Irr}(G)$, and $\theta$ a compound Clifford correspondent of $\chi$. If $\theta$ is generalized strongly monomial, then so is $\chi$.
\end{prop}
\noindent  The above proposition is  a key for the following:
\begin{theorem}\label{t3}
	Supermonomial groups are generalized strongly monomial.
\end{theorem}
\begin{theorem}\label{t4}
	If $\chi \in \operatorname{Irr}(G)$ is a  supermonomial character, then it is  generalized strongly monomial.
\end{theorem} \noindent  \noindent  Theorem \ref{t3} is a direct consequence of Theorem  \ref{t4},  which we will prove in this section. We begin by proving Proposition \ref{t2} which is also crucial for the results in the forthcoming sections.  \para \noindent{\bf Proof of Proposition \ref{t2}.} Observe that it is enough to prove the result when $\theta$ is a Clifford correspondent of $\chi$. Suppose $N \unlhd G$,  $\psi \in \operatorname{Irr}(N | \chi_N)$, and $\theta \in \I (T)$, where $T= I_{G}(\psi)$, is the Clifford correspondent of $\chi$ w.r.t $\psi$. We will show that if $\theta$ is generalized strongly monomial, then so is $\chi$. For this purpose, it is sufficient to prove the following:\begin{description} \item[(i)]  $T \unlhd  \operatorname{Cen}_{G}(e_{\Q}(\theta))$; \item[(ii)] the distinct $G$-conjugates of $e_{\Q}(\theta)$ are mutually orthogonal. \end{description} For if $\theta$ is generalized strongly monomial, then $\theta = \lambda^T$ for some linear character $\lambda$ of a subgroup $H$ of $T$ and a strong inductive chain $H=H_0 \leq H_1 \leq \cdots \leq H_m = T$ from $H$ to $T$, i.e.,   $H_i \unlhd \operatorname{Cen}_{H_{i+1}}(e_{\mathbb{Q}}(\lambda^{H_i}))$ and  the distinct $H_{i+1}$-conjugates of $e_{\mathbb{Q}}(\lambda^{H_i})$ are mutually orthogonal, for all $0 \leq i \leq m-1$. Now, if (i) and (ii) as stated above hold then we obtain that  $H=H_0 \leq H_1 \leq \cdots \leq H_m = T \leq G$ is a strong inductive chain from $H$ to $G$, showing that $\lambda^G = \chi$ is generalized strongly monomial.\par We now move to prove (i) and (ii). In the first step, we show that \linebreak $T \unlhd \operatorname{Cen}_{G}(e_{\mathbb{Q}}(\psi)).$  Let $x \in \operatorname{Cen}_{G}(e_{\mathbb{Q}}(\psi))$.  Then $e_{\mathbb{Q}}(\psi)^{x^{-1}}=e_{\mathbb{Q}}(\psi)$ which gives that $\psi^{x^{-1}}=\sigma
\circ \psi$, where $\sigma \in \operatorname{Gal}(\mathbb{Q}(\psi)/\mathbb{Q})$. Let $t \in T$. Now $x^{-1}tx \in T$ if, and only if, $\psi^{x^{-1}tx}(n) = \psi(n)$ for all $n \in N$. Note that $\psi(x^{-1}txnx^{-1}t^{-1}x)=\psi^{x^{-1}}(txnx^{-1}t^{-1})=
\sigma\circ\psi(txnx^{-1}t^{-1}) =\frac{1}{e}(
\sigma\circ\theta(xnx^{-1}))=\sigma\circ\psi(xnx^{-1})=
\psi(n),$ where $\theta_{N}= e\psi$ with $e=\langle \theta_{N}, \psi \rangle$. Therefore, it follows that  $x^{-1}tx \in T$. This proves that $T \unlhd  \operatorname{Cen}_{G}(e_{\Q}(\psi))$. \par Our next step is to prove that $\operatorname{Cen}_{G}(e_{\mathbb{Q}}(\theta)) $ is a subgroup of $\operatorname{Cen}_{G}(e_{\mathbb{Q}}(\psi))$, which will prove (i) in view of the above step. This, however, follows from the following stronger statement  which we will prove: \begin{equation}\label{e3} e_{\mathbb{Q}}(\theta)e_{\mathbb{Q}}(\theta)^{g}=0 ~~\forall~ g \in G\setminus \operatorname{Cen}_{G}(e_{\mathbb{Q}}(\psi)).\end{equation}  Due to Lemma 1 of \cite{BK2}, we have $$e_{\mathbb{Q}}(\theta)e_{\mathbb{Q}}(\psi)= e_{\mathbb{Q}}(\theta)= e_{\mathbb{Q}}(\psi)e_{\mathbb{Q}}(\theta).$$  Consequently, for any $g \in G$, \begin{equation} \label{g1} e_{\mathbb{Q}}(\theta)e_{\mathbb{Q}}(\theta)^{g}=e_{\mathbb{Q}}(\theta)e_{\mathbb{Q}}(\psi)
	e_{\mathbb{Q}}(\psi)^{g}e_{\mathbb{Q}}(\theta)^{g}. \end{equation}  As $N\unlhd G$, $e_{\mathbb{Q}}(\psi)$ and $e_{\mathbb{Q}}(\psi)^{g}$
are both primitive central idempotents of $\mathbb{Q}N$. Since $g \notin \operatorname{Cen}_{G}(
e_{\mathbb{Q}}(\psi))$, $e_{\mathbb{Q}}(\psi)$ and $e_{\mathbb{Q}}(\psi)^{g}$ can't be same and hence they are mutually orthogonal. Therefore, it follows from eqn \ref{g1} that $e_{\mathbb{Q}}(\theta)$ and $e_{\mathbb{Q}}(\theta)^{g}$ are mutually orthogonal. Hence eqn \ref{e3} is proved. \par Note that eqn \ref{e3} partially proves (ii) also. To finish the proof of (ii), it only \linebreak remains  to show  that $e_{\mathbb{Q}}(\theta)e_{\mathbb{Q}}(\theta)^{g}=0$ if $g\in \operatorname{Cen}_{G}(e_{\mathbb{Q}}(\psi))\setminus \operatorname{Cen}_{G}(e_{\mathbb{Q}}(\theta))$. Consider $g \in \operatorname{Cen}_{G}(e_{\mathbb{Q}}(\psi))$. In view of the first step, i.e., $T \unlhd \operatorname{Cen}_{G}(e_{\mathbb{Q}}(\psi)),$  we have that both $e_{\mathbb{Q}}(\theta)$ and $e_{\mathbb{Q}}(\theta)^{g}$ are primitive central idempotents of $\mathbb{Q}T$. Hence either they are same or mutually orthogonal. But if $g \not\in \operatorname{Cen}_{G}(e_{\mathbb{Q}}(\theta))$, they can't be same and therefore are mutually orthogonal. This finishes the proof of (ii) and hence completes the proof of the proposition. \qed \vspace{.5cm} \\
\noindent{\bf Proof of Theorem \ref{t4}.}  Let  $\chi \in \operatorname{Irr}(G)$  be a supermonomial character. Consider a stabilizer limit, say $\theta$, of $\chi$.  Then $\theta$ is quasi-primitive. As $G$ is solvable, by a \linebreak theorem of Berger \cite{Berger}, quasi-primitivity of a character is equivalent to the primitivity. So $\theta$ is primitive. But, the supermonomiality of  $\chi$ implies that every primitive character inducing $\chi$ must be linear. Hence $\theta$ is linear and thus generalized strongly monomial. Consequently, by Proposition \ref{t2},  $\chi$ turns out to be  generalized strongly monomial. \qed \vspace{.5cm}\\ Due to Theorem \ref{t3}, many interesting classes of groups turn out to be generalized strongly monomial:
\begin{cor}\label{c1} The following solvable groups being supermonomial are generalized strongly monomial:
	\begin{description} \item[(i)] groups whose all primitive characters are linear and every proper subgroup is monomial, in particular, subgroup closed monomial groups;
		\item[(ii)] supersolvable-by-Sylow abelian groups and Sylow abelian-by-supersolvable groups, where Sylow abelian denotes the class of solvable groups with all its Sylow subgroups abelian; \item[(iii)]  monomial groups of odd order whose irreducible characters are of prime power degree; \item [(iv)] groups whose order is a product of distinct odd primes; \item [(v)] groups of odd order with every irreducible character being a $\{p\}$-lift for possibly different primes $p$. \end{description}
\end{cor}\noindent{\bf Proof.}\noindent~ (i) Lewis, in  Lemma 2.3 of \cite{Lewis}, proved that such groups are supermonomial.\vspace{.1cm}\\(ii) Since supersolvable-by-Sylow abelian groups and Sylow abelian-by-supersolvable  groups are subgroup closed monomial, we are done, using (i). \vspace{.2cm}\\(iii) From (\cite{IsaacsV}, Theorem 9.14), such groups are supermonomial.  \vspace{.2cm}\\(iv) Let $\chi$ be an arbitrary irreducible character of a subgroup, say $H$, of $G$. \linebreak Consider a stabilizer limit, say $\theta \in \operatorname{Irr}(S)$ where $S \leq H$, of $\chi$. Then $\chi = \theta^{H}$ with $\theta$ being a primitive character. Due to Theorem 2.18 of \cite{IsaacsV}, $\theta(1)^{2}$ divides $|S:\mathcal{Z}(S)|$, where $\mathcal{Z}(S)$ denotes the center of $S$. However, the order of $G$ is a product of distinct primes. Therefore $\theta$ must be linear. This proves that every irreducible character of $H$ is monomial. Since $H$ is an arbitrary subgroup of $G$, we have that $G$ is a subgroup closed monomial group. This proves (iv), in view of (i). \vspace{.2cm}\\(v) Theorem 1 of \cite{Lewis} yields that such groups are supermonomial. \qed  \para \noindent  In response to Dornhoff's question whether arbitrary normal subgroups of monomial groups are monomial, Dade, in \cite{Dade}, provided a counterexample. In the next theorem, we will see that  this example also turns out to be generalized strongly monomial.\begin{theorem}\label{t15} Dade's example of a monomial group where monomiality is not \linebreak inherited by normal subgroups is generalized strongly monomial. \end{theorem} \noindent{\bf Proof.} The example given by Dade in \cite{Dade} is that of a group $G$  (of even order) which has a central subgroup $Z$ such that $G/Z$ is abelian-by-supersolvable. So any irreducible character of $G$ with $Z$ in its kernel is supermonomial and hence generalized strongly monomial. If an irreducible character $\chi$ of $G$ does not contain $Z$ in its kernel, then it's shown  while proving the monomiality of such character in \cite{Dade} that $\chi$ is induced from a supermonomial character, of a normal subgroup of $G$ (also see \cite{BH}, Example 24.11). It can be readily verified that if $\varphi= \vartheta^{G}$, where $\vartheta$ is an irreducible character of a normal subgroup, say $H$, of $G$ and if $\vartheta$ is generalized strongly monomial, i.e., $\vartheta = \lambda^{H}$, where $\lambda$ is a linear character of $T$, having a strong inductive chain $T=H_{0}\leq H_{1}\leq \cdots \leq H_{n}=H$ from $T$ to $H$, then $\varphi$ is also generalized strongly monomial having a strong inductive chain $T=H_{0}\leq H_{1} \leq \cdots \leq H_{n}=H \unlhd G$ from $T$ to $G$.  Consequently, $\chi$ is also generalized strongly monomial. This proves that any arbitrary irreducible character of $G$ is generalized strongly monomial.\qed
\vspace{.2cm}\\ In \cite{IsaacsII}, Isaacs proved that if $\chi$ is an irreducible complex character of an arbitrary finite solvable group $G$ such that its restriction to the fitting subgroup of $G$ is irreducible, then all the  primitive characters that induce $\chi$ must have same degree. So if, in addition, $\chi$ is monomial then it turns out to be supermonomial and thus generalized strongly monomial as well. This gives the following:
\begin{cor}\label{c2} If $\chi $  is a monomial character of a solvable group $G$ such that its restriction to the fitting subgroup of $G$ is irreducible, then it is generalized strongly monomial.
\end{cor} \noindent   Finally, we prove a slight generalization of Theorem \ref{t3}.
\begin{theorem} \label{t5} Any central extension of a cyclic group by a supermonomial group of  coprime order is generalized strongly monomial.
\end{theorem} \noindent{\bf Proof.} Let $G$ be a group with a cyclic central normal subgroup $N$ such that $G/N$ is supermonomial and $\gcd(|N|, |G/N|) =1$.  Consider $\chi \in \I (G)$. In view of  Proposition \ref{t2}, it is enough to show that there is a compound Clifford correspondent of $\chi$ which is linear. Let $\psi$ be an irreducible constituent of $\chi_{N}$.  As $N$ is central in $G$, $\psi$ is $G$-invariant. Since $\gcd(|N|, |G/N|) =1$, we have  that $\psi$ extends to $G$ (see Corollary 8.16 of \cite{IM}). Let $\delta$ be one such extension. By Gallagher's theorem (\cite{IM}, Corollary 6.17), $\chi = \delta \theta$, where $\theta \in \I (G)$ has $N$ in its kernel. Let $\overline{\theta}$ be $\theta$ modulo $N$. Since $G/N$ is supermonomial, every stabilizer limit of  $\overline{\theta}$ being quasi-primitive is primitive and hence linear. Consequently, there is a linear stabilizer limit of $\theta$. We now show that if $\theta \in \I (G)$ has a linear stabilizer limit, then so does $\delta \theta$. Let $M \unlhd G$, $\zeta \in \I (M| \theta_M)$, and $\phi\in \I (I_{G}(\zeta))$ be the Clifford correspondent of $\theta$ w.r.t $\zeta$. Let $\zeta^* = \delta_M \zeta$. Clearly $\zeta^*$ is an irreducible constituent of $(\delta\theta)_{M}.$ As $\delta$ is $G$-invariant, $I_{G}(\zeta^*)  = I_{G}(\zeta)$. Set $\phi^* = \delta_{I_{G}(\zeta)}\phi$.  One can check that $\phi^* \in \I(I_{G}(\zeta^*))$  is the Clifford correspondent of  $\delta \theta$ w.r.t. $\zeta^*$.  By the repeated application of this step, it follows that if $H$ is a subgroup of $G$ and $\lambda \in \I (H)$ is a linear stabilizer limit of $\theta$, then $ \lambda^*  = \delta_H \lambda$ is a compound Clifford correspondent of  $\delta \theta$.  Since $\lambda^*$ is linear, it has to be a stabilizer limit of $\delta \theta$ and therefore the result follows. ~\qed

\section{Groups of order $p^aq^b$} We recall some of the terminologies given by Dade and Loukaki in \cite{DL}. Let $G$ be a finite group. Consider a triple $\tau = (G,N,\psi)$  with  $N\unlhd G$ and $\psi \in \operatorname{Irr}(N)$. The \textit{center} $\mathcal{Z}(\tau)$ of $\tau$ is defined to be the center $\mathcal{Z}(\psi^{G})$ of the induced character $\psi^{G}$. A triple $\tau_{1}=(G_{1},N_{1},\psi_{1})$ is called a \textit{linear reduction} of the triple $\tau = (G,N,\psi)$ whenever there exists $L \unlhd G$ contained in $N$, and  a linear character $\lambda \in \operatorname{Irr}(L)$ lying below $\psi$, such that $G_{1}=I_{G}(\lambda)$, $N_{1}=I_{N}(\lambda)$ and $\psi_{1}$ is the unique Clifford correspondent of $\psi$ w.r.t. $\lambda$. Furthermore, a triple $\tau'$ is said to be a \textit{multilinear reduction} of $\tau$ if there is a series of triples $\tau = \tau_{0},\tau_{1},\cdots,\tau_{n}=\tau'$ such that each $\tau_{i}$ is a linear reduction of $\tau_{i-1}$ for $1\leq i\leq n$, and $\tau'$ is called a \textit{linear limit} if it is a multilinear reduction of $\tau$ in a way that the only possible linear reduction of $\tau'$ is $\tau'$ itself.\par Also, we recollect some notions introduced in \cite{CZJ} which are related to the work in \cite{DL,Loukaki2006,ML}. For a triple $\tau = (G,N,\psi)$  with  $N\unlhd G$ and $\psi \in \operatorname{Irr}(N)$, the \textit{section} of $\tau$, denoted Sec$\tau$, is defined to be the quotient group $N/\mathcal{Z}(\psi^{G})$.  It is to be noted that an arbitrary triple may have different linear limits but their sections are isomorphic to each other, in view of the main theorem of \cite{DL}. A triple $\tau$ is said to have a \textit{trivial linear limit} $\tau'$ if the section $\operatorname{Sec}\tau'$ is trivial and is said to have a \textit{nilpotent linear limit} if its section is nilpotent.  \para \noindent In the next theorem, we will show that the class of groups studied by Loukaki in \cite{ML} turns out to be generalized strongly monomial.
\begin{theorem} \label{t6} Monomial groups of order $p^a q^b$, where $p$ and $q$ are odd primes, are generalized strongly monomial. Furthermore, their normal subgroups are also \linebreak generalized strongly monomial.
\end{theorem}
\begin{lemma}\label{l1} Let $G$ be a finite group and let $\chi \in \operatorname{Irr}(G)$. If the triple $(G,G, \chi)$ has a nilpotent linear limit, then $\chi$ is generalized strongly monomial.\end{lemma}\noindent{\bf Proof.} Suppose $\chi \in \operatorname{Irr}(G)$ is such that $\tau = (G,G,\chi)$ has a nilpotent linear limit, say $\tau'= (G', G', \chi')$. Since $\chi'$ is a compound Clifford correspondent of $\chi$, in view of  Proposition \ref{t2}, it is enough to show that $\chi'$ is generalized strongly monomial. Clearly, $\mathcal{Z}(\chi')/\operatorname{ker}\chi'$ is a central subgroup of $G'/\operatorname{ker}\chi'$ (see Lemma 2.27 of \cite{IM}). The triple $\tau'$ being a nilpotent linear limit of $\tau$ yields that the factor group  $G'/ \mathcal{Z}(\chi')$ is nilpotent. Consequently, $G'/ \operatorname{ker}\chi'$ is nilpotent. Now the generalized strong monomiality of $\chi'$ follows by going modulo $\operatorname{ker}\chi'$ and using the fact that nilpotent groups being  subgroup closed monomial are also generalized strongly monomial. \qed\vspace{.5cm} \\  {\bf Proof of Theorem \ref{t6}.} In \cite{ML}, it is proved that if $G$ is a monomial group of order $p^a q^b$, where $p$ and $q$ are odd primes, then for every normal subgroup $N$ of $G$ and $\psi \in \I (N)$, the triple  $(G, N, \psi)$ has nilpotent linear limit. In particular, $(G,G, \chi)$ has a nilpotent linear limit for all $ \chi \in \I (G)$. Hence, by Lemma \ref{l1}, $G$ is generalized strongly monomial. Finally all normal subgroups, being monomial by \cite{ML},  are also generalized strongly monomial. \qed

\section{Groups with solvability length atmost 3 and fitting length atmost 2} In this section, Parks \cite{Parks} notion of successful search of the triple $(G,N,\psi)$ where $N$ is a nilpotent normal subgroup of $G$ and $\psi \in \operatorname{Irr}(N)$ plays a key role. This notion is related to that of trivial linear limits. In \cite{Parks}, Parks proved that if $G$ is a finite monomial group and $N\unlhd G$ with $N$ nilpotent and $G/N$ supersolvable of odd order then for every $\psi \in \operatorname{Irr}(N)$ there is a successful $(G,N,\psi)$ search. Observe that the existence of a successful search of $(G,N,\psi)$,  in Theorem 7.4 of \cite{Parks}, implies that it has a trivial linear limit. Following the notation of Chang, Zheng and Jin \cite{CZJ}, we write $N \propto G$ to denote those normal subgroups $N$ of  $G$ for which the triple $(G,N,\psi)$ has a trivial linear limit for each $\psi \in \operatorname{Irr}(N)$. Hence Parks's result \cite{Parks} can be rephrased as follows: if $N$ is a nilpotent normal subgroup of a monomial group $G$ such that the factor group $G/N$ is supersolvable of odd order, then $N\propto G$.
\par  Consider a monomial group $G$ having a  normal subgroup $N$ of $G$ with $N$ nilpotent such that the factor group $G/N$ is supersolvable of odd order. In this section we will show that such groups are generalized strongly monomial.
\begin{prop}\label{t7} Assume that $G$ is a finite monomial group and $N$ is a normal subgroup of $G$ with $G/N$ supersolvable of odd order. If $N\propto G$, then $G$ is a generalized strongly monomial group.\end{prop}\noindent{\bf Proof.} Let $\chi \in \operatorname{Irr}(G)$ and $\psi \in \operatorname{Irr}(N|\chi_{N})$. Let $(G', N',\psi')$  be a trivial linear limit of  $(G,N,\psi)$.  The construction of linear limits tells that there is a compound Clifford correspondent $\chi' \in \I(G')$ which lies above $\psi'$ and moreover $\psi'$ is a homogeneous restriction of $\chi'$.   By Proposition \ref{t2}, to show that  $\chi$ is generalized strongly monomial, it is enough to show that $\chi'$ is so.  As $\chi' |_{N'}$ is integral multiple of $\psi'$, $\ker\chi' \cap N' = \ker\psi'$ and so $\ker\psi' \leq \ker \chi'$. By going modulo $\ker\psi'$, it is enough to see that $\overline{\chi'}$ is a generalized strongly monomial character of $G'/ \ker\psi'$.  The triple $(G', N', \psi')$ being a trivial linear limit of $(G,N,\psi)$ yields that its section $N'/\mathcal{Z}(\psi'^{G'})$ is trivial, i.e., $N'  = \mathcal{Z}(\psi'^{G'})$. Hence $\psi'$ is linear and $N'/ \ker\psi'$ is abelian, using Lemma 2.27 of \cite{IM}. Thus, $G'/ \ker\psi'$ is abelian-by-supersolvable, due to Lemma 2.2 of \cite{CZJ}.  Since abelian-by-supersolvable groups are generalized strongly monomial, we have that $\overline{\chi'}$  is generalized strongly monomial. Consequently the same is true for $\chi'$ and subsequently for $\chi$. This finishes the proof of the theorem.\qed \vspace{.2cm}\\  As an immediate consequence, we  obtain the following:
\begin{theorem}\label{c0} Monomial groups of odd order which are nilpotent-by-supersolvable are generalized strongly monomial. In particular, monomial groups of odd order with  fitting length atmost 2 are generalized strongly monomial.\end{theorem}  \noindent Next, we will see that monomial groups of odd order which are nilpotent-by-nilpotent-by-Sylow abelian are also generalized strongly monomial. \begin{theorem}\label{t8} Monomial groups of odd order which are nilpotent-by-nilpotent-by-Sylow abelian are  generalized strongly monomial. In particular, monomial groups with solvability length atmost 3 are generalized strongly monomial.\end{theorem}
\noindent {\bf Proof.} Let $G$ be a monomial group of odd order with a nilpotent-by-nilpotent normal subgroup $N$ such that the factor group $G/N$ is Sylow abelian. Let $\chi \in \operatorname{Irr}(G)$ and let $\psi$ be an irreducible constituent of  $\chi_{N}$.  Due to Theorem C of \cite{Loukaki2006}, $(G,N,\psi)$ has a nilpotent linear limit, say $(G',N',\psi')$.  This gives  a compound Clifford correspondent $\chi'  \in \operatorname{Irr}(G')$  of $\chi$ which lies above the $G'$-invariant character $\psi' \in \operatorname{Irr}(N')$.  In view of Proposition \ref{t2}, it is enough to show that $\chi'$ is generalized strongly monomial. Observe that  $G'/N'$ is Sylow abelian.  Next $N'/ \operatorname{ker}\psi'$ is nilpotent as  $N'/\mathcal{Z}(\psi'^{G'})$ is nilpotent and $\psi'$ is $G'$-invariant.   Therefore, $G'/\operatorname{ker}\psi'$ is nilpotent-by-Sylow abelian and hence generalized strongly monomial (see Corollary \ref{c1}). Since  $\ker\psi' \leq \ker\chi' $, by going modulo $\operatorname{ker}\psi'$, it turns out that $\overline{\chi}'$ and hence $\chi'$ is generalized strongly monomial. This completes the proof.\qed
\section{Towers of groups with pairwise coprime order}
We recall some notions introduced by Chang, Zheng and Jin \cite{CZJ}. A solvable group $N$ is called an \textit{$L$-group}, if whenever  $G$ is a monomial group and $N\unlhd G$ then the triple $(G,N, \psi)$ has a nilpotent linear limit for every $\psi \in \operatorname{Irr}(N)$. Furthermore, $N$ is said to be an \textit{$\widetilde{L}$-group} if each subquotient of N is an $L$-group.
\begin{theorem}\label{t9} Let $G$ be a monomial group with a series of normal subgroups $$ 1 = G_0 \leq G_1 \leq \cdots  \leq G_{n-1} \leq G_n =G,$$ such that \begin{description} \item[(a)] $|G_i /G_{i-1}| $ are pairwise coprime for $ 1 \leq i \leq n$;
		\item[(b)] $G_i /G_{i-1}$ is an $\tilde{L}$-group for each $ 1 \leq i \leq n-1$;
		\item[(c)] all the subquotients of $G/G_{n-1}$ are also monomial.
	\end{description} Then $G$ is generalized strongly monomial and moreover all its normal subgroups are  also generalized strongly monomial.
\end{theorem}
\begin{prop}\label{p2} Assume that $G$ is a monomial group and $N \unlhd G$ with $G/N$ a subgroup closed monomial group. If $N\propto G$, and $\gcd(|N|, |G/N|) =1$, then $G$ is generalized strongly monomial.
\end{prop}
\noindent {\bf Proof.} Consider any $\chi \in \I (G)$. We have to show that $\chi$ is a generalized strongly monomial character. \par If the restriction of $\chi$ to  $N$ is the principal character, then $N$ is a subgroup of  $\ker \chi$. So by going modulo $N$, and using the fact that subgroup closed monomial groups are generalized strongly monomial (by Corollary \ref{c1}), it follows that $\overline{\chi}$ and hence $\chi$ is generalized strongly monomial   \par Assume that the restriction of $\chi $ to  $N$ is not the principal character. Let $ \psi$ be an irreducible constituent of $\chi_{N}$.  Since $N\propto G$, the triple $(G, N,  \psi)$ has a trivial linear limit, say $(G', N', \psi')$.  Suppose $\chi' \in \operatorname{Irr} (G')$ which lies above $\psi'$ is a compound Clifford correspondent of $\chi$. To prove that $\chi$ is generalized strongly monomial,  it is enough to show that $\chi'$ is generalized strongly monomial. Observe that  $\ker \psi' \leq \ker \chi'$. As $N'=\mathcal{Z}(\psi')$ and $\psi'$ is $G'$-invariant, we see that $G'/ \ker\psi'$  has a cyclic central subgroup $N'/ \ker\psi'$, using Lemma 2.27 of \cite{IM}. Also, its quotient is subgroup closed monomial. Furthermore, the order of $N'/ \ker\psi'$  is coprime to its index in $G'/ \ker\psi'$. Hence  Theorem \ref{t5} is applicable and  we get that $\chi'$ modulo $\ker \psi' $ is generalized strongly monomial and consequently so is $\chi'$.  This finishes the proof. \qed
\begin{cor} \label{c11}Any extension of an abelian group by a subgroup closed monomial group of  coprime order is  generalized strongly monomial.\end{cor}\noindent {\bf  Proof of Theorem \ref{t9}.} By Theorem A of \cite{CZJ}, $G_{n-1}  \propto G$ and so Proposition \ref{p2} gives that $G$ is generalized strongly monomial. Furthermore, by Theorem B of \cite{CZJ}, every normal subgroup of $G$ is monomial and thus the hypothesis of the theorem is preserved by its normal subgroups and hence they are generalized strongly monomial as well. \qed
\begin{cor}\label{c3}  Let $G$ be a monomial group with a series of normal subgroups $$ \{e\} = G_0 \leq G_1 \leq \cdots  \leq G_{n-1} \leq G_n =G,$$ such that \begin{description} \item[(a)] $|G_i /G_{i-1}| $ are pairwise coprime for $ 1 \leq i \leq n$;
		\item[(b)] $G_i /G_{i-1}$ is either nilpotent(odd~order)-by-nilpotent or Sylow abelian-by-nilpotent for  $ 1 \leq i \leq n-1$;
		\item[(c)] all the subquotients of $G/G_{n-1}$ are also monomial.
	\end{description} Then  all  normal subgroups  of $G$ are  generalized strongly monomial.
\end{cor}\noindent  {\bf Proof.}  Lemma 3.4 of \cite{CZJ} yields that nilpotent(odd~order)-by-nilpotent or Sylow abelian-by-nilpotent are $\tilde{L}$ groups. Hence the result immediately follows from the above theorem.  \qed \vspace{.5cm} \\ Another immediate consequence of Theorem \ref{t9} is the following:
\begin{cor}\label{c4} If $G$ is a monomial group with a Sylow tower, then all its subnormal subgroups are  generalized strongly monomial.
\end{cor} \noindent {\bf Proof.} Clearly a monomial group $G$ with a Sylow tower is of the type discussed in Theorem \ref{t9} and hence it is generalized strongly monomial. Also, by the work of Gunter \cite{Gunter}, every subnormal subgroup of $G$ is monomial  and thus being with a Sylow tower it is also generalized strongly monomial. \qed
 \begin{remark}\label{r1} In \cite{Fukushima}, Fukushima gave an  example  of a monomial group which has a Hall subgroup that is not monomial. Since the example of Fukushima is that of a monomial group with a Sylow tower, it turns out to be generalized strongly monomial as well. \end{remark}
\begin{theorem}\label{t10} Let $G$ be a monomial group with a series of normal subgroups \linebreak  $1 \leq L \leq K \leq G$ such that: \begin{description} \item[(a)] the factor group $G/K$ is supersolvable, $K/L$ is nilpotent and $L$ has a normal series $1=L_0 \leq L_1 \leq \cdots \leq L_{n-1} \leq L_n =L$ where all factors $L_i/L_{i-1}$ are abelian with pairwise coprime orders, for $i=1, \cdots, n$;
		\item[(b)] either $|G/K|$ or $|K/L|$ is odd;
		\item[(c)] $|K/L|$ and $|L|$ are coprime.
	\end{description} Then all normal subgroups of $G$ are  generalized strongly monomial.   Furthermore, if $L$ is abelian then all Hall subgroups of $G$ are also generalized strongly monomial. \end{theorem}
\noindent  {\bf Proof.} Consider a group $G$ which have a series of normal subgroups as mentioned in the statement and $\chi \in \I (G)$. If $K \leq \ker\chi$, then by going modulo $K$ and using the fact that supersolvable groups are strongly monomial (and hence generalized strongly monomial) it follows that $\chi$ is so. If $K \nleq \ker\chi$, consider the triple $(G, K, \psi)$, where $\psi$ is an irreducible constituent of $\chi_K$. By Theorem A of \cite{ZJ}, $(G, K, \psi)$ has a linear limit, say $(G', N', \psi'),$ where $\psi'$ is linear. There now exists $\chi' \in \I(G'|\psi') $ which induced to $G$ is $\chi$ and is its compound Clifford correspondent. The generalized strong monomality of $\chi'$ is seen by going modulo $\ker\psi'$ and using the fact that cyclic-by-supersolvable are again supersolvable and hence generalized strongly monomial. Finally, the generalized strong monomiality of $\chi$ follows from that of $\chi'$. \par Next, by Theorem B of \cite{ZJ}, all normal subgroups of $G$ are monomial and consequently being of the same type are also generalized strongly monomial. \par Finally, if $L$ is abelian then, by Theorem C of \cite{ZJ}, all Hall subgroups of $G$ are monomial and thus the result follows. \qed  \section{Questions} Dornhoff's question regarding whether arbitrary normal subgroups of monomial groups are monomial was answered in negative by Dade \cite{Dade} by constructing an example of a monomial group having a non monomial normal subgroup of even order and having even index in the group. The prime 2 played a key role in Dade's example and thus the question regarding monomiality of normal subgroups of odd order monomial groups was left open. This question has been studied extensively in literature. For various classes of monomial groups, an affirmative answer to this question has been provided. In the earlier sections, we have seen that all such monomial groups where monomiality is inherited by normal subgroups turn out to be generalized strongly monomial. This leads one to ask the following:
\begin{quote} {Question 1:  Is a normal subgroup of a generalized strongly monomial group of odd order itself monomial?} \end{quote} As mentioned earlier, Isaacs conjectured that every
monomial group of odd order is a supermonomial group. In light of Theorem \ref{t3}, where we have shown that supermonomial groups  are generalized strongly monomial, a weak form of Isaacs's conjecture is raised:
\begin{quote} {Question 2: Whether every generalized strongly monomial group of odd order is supermonomial? }  \end{quote} Another question which is of interest is  the following:\begin{quote} {Question 3: Is an arbitrary monomial group of odd order  generalized strongly monomial ? }  \end{quote} An affirmative answer to both Questions 2 and 3 is equivalent to proving Isaacs conjecture and it also answers the open question regarding monomiality of normal subgroups of odd order monomial groups in affirmative,  in view of Lewis's work in \cite{Lewis}.\bibliographystyle{amsplain}
\bibliography{BK4}

\end{document}